\documentclass[final,3p,times,authoryear]{elsarticle}

\usepackage{amssymb}
\usepackage[dvipsnames]{xcolor}

\journal{Safety Emergency Science}

\begin{document}

\begin{frontmatter}



\title{Recent Advances in Disaster Emergency Response Planning: Integrating Optimization, Machine Learning, and Simulation}


\author[inst1]{Fan Pu}
\author[inst1]{Zihao Li}
\author[inst2]{Yifan Wu\corref{cor1}}

\author[inst1]{Chaolun Ma}

\author[inst2]{Ruonan Zhao}

\affiliation[inst1]{organization={Department of Civil and Environmental Engineering, Texas A\&M University},
            city={College Station},
            postcode={77843}, 
            state={TX},
            country={United States}}

\affiliation[inst2]{organization={Department of Industrial \& Systems Engineering, Texas A\&M University},
            city={College Station},
            postcode={77843}, 
            state={TX},
            country={United States}}

\cortext[cor1]{Corresponding author. Email: wuyifan1995@tamu.edu.}

\begin{abstract}
The increasing frequency and severity of natural disasters underscore the critical importance of effective disaster emergency response planning to minimize human and economic losses. This survey provides a comprehensive review of recent advancements (2019–2024) in five essential areas of disaster emergency response planning: evacuation, facility location, casualty transport, search and rescue, and relief distribution. Research in these areas is systematically categorized based on methodologies, including optimization models, machine learning, and simulation, with a focus on their individual strengths and synergies. A notable contribution of this work is its examination of the interplay between machine learning, simulation, and optimization frameworks, highlighting how these approaches can address the dynamic, uncertain, and complex nature of disaster scenarios. By identifying key research trends and challenges, this study offers valuable insights to improve the effectiveness and resilience of emergency response strategies in future disaster planning efforts.
\end{abstract}


\end{frontmatter}

\section{Introduction}

Over the past few years, there has been an increase in both the occurrence and intensity of natural disasters like hurricanes, earthquakes, and tsunamis, leading to considerable human and financial damages \citep{douris2021atlas,li2023prediction,zhou2021analyzing}. By November 1, 2024, the United States had encountered 24 distinct weather and climate disaster events in 2024, each inflicting financial losses over \$1 billion and collectively resulting in 418 deaths \citep{NCEI}. Hurricanes Helene and Milton ravaged the southeastern United States in September and October, with the former causing approximately \$250 billion and the latter \$50 billion in damages. Combined, these hurricanes resulted in no fewer than 250 fatalities \citep{wilcox2024hurricane}. The Smokehouse Creek Fire burned an estimated 1,058,482 acres across Texas and Oklahoma, making it the second-largest wildfire in U.S. history since 1988 and prompting a state of disaster declaration in 60 counties \citep{yoon2024wildfires}. Hurricanes and wildfires are just a glimpse of the many natural disasters we face. The increasing frequency and severity of disasters highlight the critical importance of emergency response planning across five key areas, such as \textbf{evacuation}, \textbf{facility location}, \textbf{casualty transport}, \textbf{search and rescue}, and \textbf{relief distribution}. Effective planning in these areas saves lives by enabling swift and organized evacuations, optimizing the placement of emergency facilities, and ensuring timely transport of casualties to medical care. It also reduces losses by improving the efficiency of search and rescue operations and ensuring equitable and rapid relief distribution. By addressing these critical components, comprehensive emergency response planning enhances disaster preparedness, strengthens community resilience, and accelerates recovery efforts.

This paper reviews recent advancements in the above-mentioned five critical topics within disaster emergency response planning. These topics collectively encompass both pre-disaster and post-disaster operations. Pre-disaster operations involve activities undertaken before the occurrence of a disaster, focusing on preparation and mitigation strategies to minimize potential damage and ensure readiness. Post-disaster operations address the immediate and short-term responses following a disaster, aiming to alleviate its impacts and provide urgent assistance to affected populations.

  Following previous review papers in this field \citep{altay2006or, caunhye2012optimization, galindo2013review}, our paper emphasizes operational strategies based on optimization models such as mixed-integer programming (MIP), stochastic programming, and robust programming. Optimization models use mathematical formulations to find the best possible solution to a problem. One key strength of the optimization models is that they can provide optimal or near-optimal solutions for problems with clear objectives and constraints. However, these optimization models face challenges in addressing large-scale problems, considering uncertainty and complex environments inherent in disaster scenarios. In particular, precise problem definitions may not always be feasible in disaster scenarios.


To address these limitations, we examine the growing potential of machine learning (e.g., supervised learning, reinforcement learning (RL)) and simulation technologies (e.g., agent-based model (ABM), discrete-event simulation) in disaster emergency response planning as complementary or alternative approaches. 
Machine learning can learn patterns from data to make predictions or decisions, making it particularly effective in disaster scenarios with incomplete or uncertain information. Simulation involves creating models that replicate real-world systems, enabling the testing of various scenarios and forecasting outcomes. It has the advantage of managing the dynamic nature of disaster environments. While machine learning and simulation effectively capture the complexities of disasters, they exhibit limitations in providing optimal operations.
As a result, our key contribution lies in exploring the role of machine learning and simulation technologies alongside optimization models in disaster emergency planning. By comparing the focus of these approaches, we aim to comprehensively understand the challenges associated with disaster response, ultimately offering valuable insights into improving emergency planning strategies.

This paper reviews literature published between 2019 and 2024, encompassing the five key topics mentioned earlier, as shown in Figure 1. For each topic, we categorize studies based on their methodologies (optimization model, machine learning, and simulation), if possible, and provide a detailed review of the most recent advancements. Publications are searched from reputable publishers, including INFORMS, Elsevier, Springer, IEEE, and databases such as Web of Science and Google Scholar. It is important to note that this survey does not aim to be exhaustive; instead, it seeks to highlight and synthesize some of the most representative studies in the field. Additionally, we excluded research focused on small-scale emergency planning scenarios, such as indoor, building, or community-level cases. Papers unrelated to the operational strategy in machine learning and simulation (e.g., image recognition, device communication, and signal processing, etc.) are also omitted. In the following section, we will break down the five key areas in detail, with each section expanding on the main approaches, including optimization models, machine learning, and simulation. For areas where no relevant studies are available, they will be excluded accordingly, such as machine learning in facility location and simulation in search and rescue. The list of abbreviations used in this paper is presented in Table \ref{t:abbre}.


\begin{figure}[h]
\centering
\includegraphics[width=0.8\textwidth]{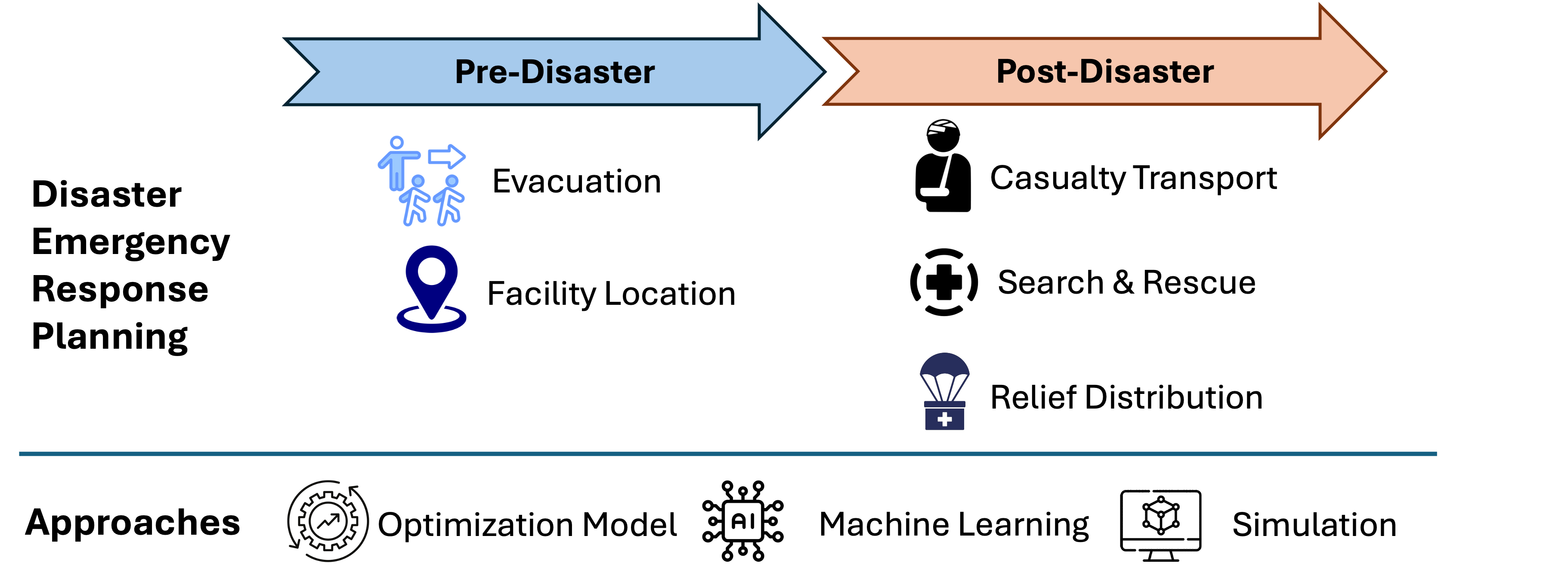}
\caption{Review Framework of Disaster Emergency Response Planning}
\label{fig:vs}
\end{figure}


\begin{table}[h]
\centering
\caption{Abbreviations and definitions.}
\label{t:abbre}
\begin{tabular}{cl}
\hline
Abbreviation & Definition                 \\ \hline
ABM                  & Agent-Based Model                   \\
CTP                  & Casualty Transport Problem           \\
GIS                  & Geographic Information System       \\
MCI                  & Mass Casualty Incident              \\
MEDEVAC              & Military Medical Evacuation Problem \\
MDP                  & Markov Decision Process             \\
MIP                  & Mixed-Integer Programming           \\
OP                   & Orienteering Problem                \\
RL                   & Reinforcement Learning              \\
UAV                  & Unmanned Aerial Vehicle             \\
VaR                  & Value-at-Risk                       \\
VRP                  & Vehicle Routing Problem             \\
FLP                  & Facility Location Problem            \\
SAR                  & Search and Rescue                    \\\hline
\end{tabular}
\end{table}


\section{Review of Evacuation}\label{LB:EVA}
Disaster evacuation is a key focus in pre-disaster operations, aimed at minimizing human casualties and mitigating property damage during natural disasters. Despite extensive research, disaster evacuation remains a complex challenge due to several factors.
First, large-scale evacuations executed without adequate strategies can overwhelm transportation systems, resulting in severe traffic congestion and delays in securing safe passage. Second, the inherent unpredictability of disasters introduces uncertainty and dynamic conditions, requiring both proactive planning and real-time adjustments. Furthermore, behavioral and psychological factors play a critical role in evacuation scenarios. Variations in risk perception among individuals lead to differing participation levels, evacuation timings, and route preferences, which further complicate evacuation planning. Finally, evacuation operations must address the needs of vulnerable populations, including the elderly, individuals with disabilities, and those without access to private transportation. These groups face disproportionate challenges, necessitating targeted resources and accommodations.
This review synthesizes the existing literature addressing these challenges, highlighting strategies and frameworks proposed to improve disaster evacuation effectiveness.

\subsection{Optimization Model in Evacuation}


The optimization of evacuation strategies has been extensively studied in recent literature. \citet{hafiz_hasan_largescale_2021} provided a systematic examination of the large-scale zone-based evacuation planning problem, where the evacuation area is divided into zones to enhance public communication and minimize confusion. Variants of the zone-based evacuation problem, incorporating contraflow, convergent plans, and non-preemption, are modeled using an MIP framework. These models are addressed using Benders decomposition, which divides the original problem into a master problem and multiple smaller subproblems of manageable size \citep{yin2023integrated,pu2022rolling}, as well as a column generation approach. Their performance is further evaluated through simulations detailed in a companion paper \citep{hasan_span_2021}. Additionally, \citet{li_tsunami-induced_2019} proposed a hierarchical evacuation structure to manage vehicle routing during tsunami evacuations. This study employs lexicographic minimax optimization combined with Tabu search techniques to derive optimal solutions.

Some studies focus on addressing uncertainties during the evacuation process. \citet{karabuk_multi-stage_2019} introduced a multi-stage stochastic optimization framework to manage evacuations under spatial and temporal weather uncertainties. A tornado forecast framework provides staged projections of speed, intensity, and direction to guide evacuation strategies. Traffic flow dynamics are modeled using the cell transmission model. Similarly, \citet{lim_robust_2019} addressed uncertainty in evacuation demand using a robust chance-constrained model. Instead of relying on exact demand distributions, they utilize partial information about the demand distribution. A path-based model is developed to determine the minimum clearance time. \citet{darvishan_dynamic_2021} proposed a real-time rerouting strategy for evacuation networks experiencing road disruptions. A dynamic network flow optimization model is formulated to minimize evacuation clearance time by adapting to real-time conditions.

Other studies concentrate on evacuating populations with special needs. \citet{rambha_stochastic_2021} address hospital evacuations during hurricanes. They develop a stochastic programming model based on a scenario tree that considers risks, transportation costs, and staffing expenditures. The model determines the evacuation order of patients over time, updating decisions as the hurricane evolves. \citet{moug_shared-mobility-based_2023} focused on populations without access to private vehicles. They propose a ride-sharing-based evacuation model, recruiting volunteer drivers to transport evacuees. Similarly, \citet{zhao_round-trip_2020} developed a round-trip bus evacuation model to optimize bus scheduling and routing for evacuating car-free individuals. This model effectively coordinates pick-up and drop-off locations to ensure efficient transportation.

\subsection{Machine Learning in Evacuation}

Machine learning offers powerful tools for predicting behavioral factors and traffic conditions by utilizing historical data. These predictions can significantly enhance evacuation strategies and preparedness. Several studies focus on behavioral analysis using machine learning. Based on survey datasets collected from hurricane evacuees, \citet{chen_data-driven_2024} examined evacuation compliance under government orders, while \citet{sun_predicting_2024} and \citet{karampotsis_understanding_2024} investigated evacuation decision-making processes. \citet{sreejith_modelling_2022} analyzed evacuation preparation time to understand how individuals respond to evacuation scenarios. Additionally, \citet{rahman_deep_2023} developed a deep learning framework that learns spatiotemporal traffic patterns at the network scale. Their model is trained on non-evacuation traffic data, offering insights into traffic behavior during evacuation scenarios.

Some studies integrate machine learning with optimization models to improve decision-making. \citet{yang_machine_2021} propose a machine learning-based flight dispatch system for large-scale evacuations. The problem is formulated as a weighted bipartite graph-matching model, where RL is employed to estimate edge weights based on historical flight data. \citet{bagloee_hybrid_2019} address contraflow optimization using a bilevel optimization framework. A hybrid heuristic method that combines machine learning with optimization techniques is developed. The machine learning model iteratively improves contraflow scenarios, enhancing the decision-making process.

\subsection{Simulation in Evacuation}

Optimization models often simplify the representation of traffic dynamics and human behavior to improve the tractability of evacuation planning. However, simulations provide a more detailed and realistic depiction of evacuation processes, capturing the complexity of human interactions and traffic dynamics. ABM is a widely used computational approach that simulates the actions and interactions of individual evacuees, represented as agents, during evacuation scenarios.

Several studies have adopted ABMs to analyze evacuation dynamics in specific contexts. \citet{fathianpour_tsunami_2023}, \citet{kim_agent-based_2022}, and \citet{takabatake_tsunami_2020} employ ABMs to evaluate tsunami evacuations, focusing on pedestrian movements and the impact of evacuation routes. \citet{keykhaei_multi-agent-based_2024} propose a situation-aware evacuation framework using ABMs to model the behavior of pedestrians, vehicles, and controlling agents. Their work examines both the interactions between agents, such as the coordination between pedestrians and vehicles, and the interactions between agents and their environment, including obstacles and evacuation signals. Similarly, \citet{harris_agent-based_2022} introduced FLEE, a modeling framework based on ABMs, to study the interconnected dynamics of hurricane forecast-warning systems. This framework integrates models of natural hazards, human systems, built environments, and their interactions, providing a comprehensive tool for analyzing the effectiveness of warning systems and evacuation decisions. In addition, \citet{na_agent-based_2019} developed an agent-based discrete-event simulation framework equipped with an embedded Geographic Information System (GIS) module for no-notice natural disaster evacuation planning. This framework incorporates multiple types of evacuees, including individuals with varying levels of mobility and access to vehicles, and models interactions between different types of vehicles, such as buses and private cars, to enhance the accuracy of evacuation planning.

Beyond addressing ``what if'' scenarios, some studies combine simulation with optimization methods to address the complexities of evacuation planning. \citet{idoudi_agent-based_2022} propose a framework that integrates the shelter allocation problem with the dynamic traffic assignment problem of evacuation. A dynamic ABM simulator captures traffic dynamics during evacuation, while the dynamic traffic assignment problem is formulated under the stochastic user equilibrium principle to account for individual route choices in uncertain conditions. This integrated approach ensures the simultaneous optimization of shelter assignments and traffic flow. Similarly, \citet{alam_evacuation_2022} developed a framework combining optimization and simulation to plan mass evacuations for individuals with mobility needs. The framework uses a Monte Carlo simulation method to estimate the time distribution of emergency vehicle dwell times at pick-up locations. A MIP model is then applied to optimize the allocation of emergency vehicles and evacuation destinations. Finally, a traffic evacuation micro-simulation model evaluates the impact of these strategies through scenario analysis, providing insights into evacuation efficiency and resource allocation.


In summary, most optimization models in evacuation research focus on developing new evacuation strategies, addressing the uncertainty of disasters, and considering populations with special needs. Stochastic programming and robust optimization are commonly used to manage uncertainties, including disaster dynamics, evacuation demands, and road disruptions. To ensure the feasibility of modeling, traffic conditions and human behavior are often simplified in these optimization frameworks, which may limit the models’ realism. Machine learning and simulation techniques help bridge these gaps. The two methods leverage historical data and complex evacuation systems for more realistic evacuation performance. However, most machine learning and simulation research reviewed are primarily used as evaluation tools and lack systematic guidance for evacuation decision-making. Consequently, there is a growing interest in integrating these methodologies to enhance the overall effectiveness of evacuation planning and decision support.

\section{Review of Facility Location}

In disaster emergency response planning, the facility location problem (FLP) aims to identify optimal strategic locations for facilities such as warehouses, shelters, medical centers, and distribution hubs to facilitate other emergency response operations such as relief distribution, casualty transport, and evacuation. Consequently, studies focusing solely on the FLP are relatively uncommon. Instead, the FLP is often integrated with other operational strategies. Such integration arises from the fact that subsequent operational strategies heavily depend on the effectiveness of location decisions, emphasizing the importance of FLP.

\subsection{Optimization Model in Facility Location}

\citet{sanci2021integer} study the integrated facility location and network restoration problem, addressing the coordination of relief distribution and network restoration efforts. They propose an integer L-shaped algorithm, utilizing branch-and-cut procedures to determine optimal facility locations and restoration resource allocations efficiently. \citet{liu2019integrated} developed a bi-objective optimization model for temporary medical facility location and casualty allocation, focusing on maximizing expected survival rates and minimizing operational costs. \citet{monemi2021multi} addressed the multi-period hub location problem for humanitarian aid distribution in refugee camps. They propose a hub-and-spoke network design where hubs serve as distribution centers to manage aid delivery efficiently across multiple periods. A compact mathematical model, enhanced with Benders decomposition, is tested with real and synthetic datasets, showcasing its effectiveness in optimizing aid logistics.

 Some papers focus on handling uncertainty. A widely used approach is robust optimization, which considers worst-case scenarios to ensure reliable solutions under all conditions within a defined uncertainty set. For instance, \citet{zhang2023improving} combined the shelter location problem with the capacity allocation problem using a scenario-wise distributionally robust optimization approach. The model integrates uncertainties by leveraging correlations between disaster severity and casualty numbers through a scenario-based ambiguity set. \cite{yin2024two} proposed a two-stage recoverable robust optimization model to optimize worst-case scenarios for integrated facility location, resource allocation, and evacuation planning. The model incorporates uncertainties in disaster impacts using a discrete scenario set that captures variations in facility capacities, relief demands, and transportation link availability. 
 
 Other studies manage uncertainties through Value-at-Risk (VaR), quantifying potential losses under a specified confidence level. \citet{zhong2020risk} incorporated conditional VaR with regret into their bi-objective model for disaster relief facility location integrated with vehicle routing. This approach measures the expected regret in worst-case scenarios and overall system performance, balancing cost and waiting time while addressing the decision-maker's risk aversion. Similarly, \citet{jin2024risk} employed worst-case mean conditional VaR in a drone-supported relief FLP. This risk-averse distributionally robust optimization model manages uncertain demands by adopting ambiguity sets to describe incomplete demand distributions.

Beyond uncertainty considerations, equity and human-centric factors are increasingly pivotal in facility location optimization for disaster emergency response. These considerations prevent resource allocation from disadvantaging undeserved or vulnerable communities. \cite{esposito2021integrated} addressed the shelter location and evacuation routing problem with a two-stage MIP model. The model integrates self- and supported-evacuation modes while considering road disruptions and equity in resource allocation. Equity is defined as ensuring fair access to shelters within a maximum travel time. Similarly, \cite{seraji2022integrative} tackled distributive injustice in humanitarian logistics using a two-stage location-allocation model to optimize pre-disaster shelter and depot locations. Their equity criterion minimizes disparities in resource availability by prioritizing sites that reduce access inequities while maintaining efficiency.

\subsection{Simulation in Facility Location} 

The simulation addresses facility location challenges in disaster emergency response by dynamically managing complexities and uncertainties. It enables real-time demand estimation, adapts to changing conditions, and integrates infrastructure dynamics such as transportation and communication systems. 


For tsunami scenarios, \cite{mostafizi2019agent} highlighted the importance of simulation in vertical evacuation shelter location planning. Their ABM framework incorporates evacuee behaviors, such as word-of-mouth communication and impatience, alongside tsunami wave dynamics to assess how shelter proximity to population centers impacts mobility and survival rates. Results demonstrate that strategically placed shelters increase participation in vertical evacuations, significantly reduce mortality, and expand viable site options. Conversely, poorly located shelters heighten risks by prolonging evacuation times and complicating decision-making for evacuees. In cascading disaster scenarios, \cite{khalili2022stochastic} underscored the potential of simulation to enhance decision-making under uncertainty. They develop a two-echelon bi-level stochastic model combining simulation with GIS-based analysis to optimize distribution center locations. The model dynamically estimates the demand for critical supplies like water, food, and medicine by analyzing interactions among vulnerable infrastructure systems and cascading effects. This approach prioritizes critical areas, informed resource allocation, and supports both pre-disaster preparedness and post-disaster response, showcasing the adaptability and precision of simulation for disaster logistics.


In brief, optimization models for the FLP are often integrated with other operational strategies, such as relief distribution and evacuation planning. A significant portion of the research focuses on disaster risk management, employing methodologies such as robust optimization to address worst-case scenarios and VaR to quantify potential losses under specified confidence levels. Additionally, many models incorporate equity criteria to minimize disparities in resource allocation, ensuring fair access to critical services. On the other hand, simulation enhances location optimization efforts by capturing the dynamic and real-time aspects of disaster response, such as evacuee behaviors, which are often oversimplified in traditional optimization models. However, integrating simulation into optimization frameworks often sacrifices solution quality, as it may result in poor local optimal solutions.

\section{Review of Casualty Transport} \label{LB:Casualty}

The casualty transport problem (CTP) focuses on efficiently transferring injured individuals from incident sites to emergency medical centers (hospitals and alternative care facilities) during mass casualty incidents (MCIs). It plays an important role in minimizing treatment delays, providing services to casualties, and preventing loss of life after disasters. Due to the suddenness and severity of disasters, CTP faces significant challenges, including time constraints, resource shortages, and scheduling difficulties. The primary focus of CTP is to respond to as many injured individuals as possible in a timely manner, using constrained resources efficiently. A closely related topic, known as the military medical evacuation problem (MEDEVAC), addresses the evacuation and dispatch of injured soldiers on the battlefield. While MEDEVAC is tailored explicitly for military applications, it shares common methodologies and concepts with disaster-related CTP. Therefore, our survey incorporates relevant MEDEVAC literature to provide a more comprehensive overview of these interconnected fields.

\subsection{Optimization Model in Casualty Transport}

Many studies have explored triage and casualty prioritization in the context of CTP \citep{feng_casualty_2021-1}. \citet{shi_treatment_2024} examined a treatment scheduling model for MCIs, incorporating victim health deterioration and wait-dependent service times. Their study identifies the conditions under which treatment priority should be given to delayed patients, offering insights into the interplay between urgency and timing. Similarly, \citet{shin_emergency_2020} investigated the integration of patient prioritization and hospital selection to maximize the expected number of survivors in MCIs. Using a Markov decision process (MDP) model, the study determines which patient should be transported first and to which hospital, leveraging real-time hospital state information. Additionally, a heuristic policy is proposed to facilitate practical implementation in emergency scenarios.

The transportation of casualties has also garnered significant research attention. \citet{chou_emergency_2022-1} develop a non-linear programming model to minimize transfer delays caused by traffic congestion. Their model considers both ambulance routing and the operational status of hospitals in response to MCIs, utilizing the cell transmission model to represent traffic dynamics. In a related study, \citet{pei_dynamic_2023} addressed the impact of potential travel delays within road networks, further emphasizing the importance of accounting for dynamic traffic conditions in casualty transport. Several studies also focus on ambulance dispatch. For instance, \citet{yoon_dynamic_2021-1} and \citet{dubois_dispatching_2021} formulate MDP models to optimize ambulance dispatch decisions, ensuring that prioritized patients receive timely medical attention.

Uncertainty in casualty transport is another critical area of research. \citet{sun_novel_2022} address risks associated with disruptions in temporary medical centers and the uncertain number of casualties in various scenarios. They develop a scenario-based robust optimization model to jointly plan facility locations, allocate relief resources, and manage casualty transport under uncertainty. In a different approach, \citet{shiri_online_2024} proposed an online optimization framework for ambulance routing in road networks. This model is specifically designed to handle uncertainties in travel times, triage levels, and treatment durations, offering a flexible and adaptive solution to real-time challenges.

\subsection{Machine Learning in Casualty Transport}

Machine learning has been extensively applied to casualty triage and dispatch problems \citep{xia_emergency_nodate, lee_multi-agent_2021-1}, for these challenges can often be modeled as optimal control problems with MDP model. However, the increasing complexity of CTP makes it extremely difficult to derive analytical solutions in the MDP framework. Consequently, RL algorithms have gained recognition as a promising approximation method for solving MDP-based models \citep{rodriguez_solving_2023, lee_improving_2020-1, robbins_approximate_2020, ma2022lyapunov,allison2024performance}. For instance, \citet{jenkins_approximate_2020} utilized approximate dynamic programming methods to solve the MEDEVAC problem. Two learning techniques based on least-squares temporal difference and a neural network are compared for policy evaluation. Results show that these learning-based technologies perform better than the traditional policy benchmark. \citet{jenkins_approximate_2021} later extended the approximate dynamic programming research to solve the combined MEDEVAC dispatching, preemption-rerouting, and redeployment problem. These techniques enhance the adaptability of casualty transport systems by learning from real-time data and optimizing outcomes in uncertain environments.


\subsection{Simulation in Casualty Transport}

Simulation models provide a platform to replicate the complex dynamics of casualty transport systems. They enable performance evaluation, strategy analysis, and resource allocation under realistic conditions. By modeling interactions among patients, resources, and infrastructure, simulations offer actionable insights to improve casualty transport processes and support decision-making in disaster scenarios. \citet{chang_simulation_2023} introduced a hybrid simulation-optimization approach, combining the discrete event simulation model with a two-stage sequential algorithm to optimize the casualty collection point location and resource allocation problem. Simulation plays an important role in handling the dynamic nature of CTP, including time-varying casualty arrivals, random triage service times, and stochastic travel times caused by road network vulnerability. \citet{jat_mass-casualty_2020} use discrete event simulation to examine the casualty-load distribution in MCIs, revealing the trade-off between travel time and queuing delays. This study criticizes the inefficiency of policies focused solely on hospital capacity.

To sum up, optimization models in casualty transport focus on triage prioritization, ambulance dispatch, and routing under complex constraints, such as time, resources, and traffic conditions. Among these, MDPs are commonly used as the modeling framework for casualty transport. However, as the problem scale increases, traditional optimization models become computationally infeasible. In such cases, RL has emerged as a powerful tool to address the growing complexity of CTP. Despite challenges related to interpretability, RL can often converge to an effective policy when sufficient training data is available. Simulation models replicate the complex casualty transport system to evaluate performance, analyze strategies, and provide insights for improving casualty transport processes. However, the optimization component is often weakened when the simulation is combined with optimization models.

\section{Review of Search and Rescue}\label{LB:SAR}

Search and Rescue (SAR) is a critical post-disaster operation focused on locating, accessing, and recovering individuals trapped, injured, or missing due to natural catastrophes. Given the need for immediate response, SAR operations often face significant resource limitations in the early stages. For example, the survival window is typically limited to 72 hours, and the availability of SAR teams is restricted, so planning the search areas and allocating search time are crucial considerations. Additionally, SAR is highly time-sensitive, as rescuing injured individuals is a top priority, with survival rates declining rapidly over time. Furthermore, the uncertainty of post-disaster environments makes rapid responses to unforeseen conditions necessary. These challenges have driven continuous research in this field.

It is worth noting that the operational planning problem of SAR is often modeled as an orienteering problem (OP). The OP involves determining optimal routes for vehicles with limited resources, such as time, to maximize the total value of visited locations. In addition, unmanned aerial vehicles (UAVs) or drones have become widely utilized in post-disaster SAR due to their high mobility, ability to gather aerial information, and access to hazardous areas. Consequently, besides SAR problems, this survey also includes studies related to the OP and the UAV-based search and routing problems. These advancements provide valuable insights that can be directly applied to SAR planning.

\subsection{Optimization Model in SAR}

We begin by reviewing the UAV-based search and routing problem. \citet{moskal2019adaptive} addressed the UAV information collection problem by modeling it as a prize-collecting vertex routing problem. This study provides a foundation for optimizing UAV missions by focusing on maximizing the value of collected information. Extending this work, \citet{moskal_unmanned_2023} investigated a more complex scenario involving regions with stochastic attributes. Uncertainties in this problem arise from factors such as information availability, sensor effectiveness, and travel times between regions. They propose a MIP model to maximize expected information collection while mitigating the risk of mission delays, offering a robust framework for UAV operations under uncertainty. Similarly, \citet{glock_mission_2020} focused on mission planning for UAV routing to sampling locations in post-incident scenarios. Their approach predicts the distribution of hazardous substances across affected areas, highlighting the potential of UAVs in environmental monitoring during disasters.

\citet{bravo_use_2019} tackled the UAV search problem using a partially observable MDP. This framework is particularly effective for disaster scenarios where imperfect information is common. The model incorporates a mechanism to update environmental understanding as UAVs collect new data, enhancing decision-making during search operations. \citet{poikonen_mothership_2020} introduced the mothership-and-drone routing problem, focusing on a two-vehicle tandem where a large mothership navigates continuous Euclidean space while carrying drones. They propose both an exact branch-and-bound approach and a heuristic method to solve this problem for large instances. Similarly, \citet{morandi_orienteering_2024} extended the classical OP to incorporate multiple drones working in coordination with a truck. Drones, constrained by limited battery life, can either travel alongside the truck without energy consumption or be deployed for short flights. A MIP formulation and a tailored branch-and-cut algorithm are developed to optimize this cooperative routing problem.

One notable advancement in OP is the consideration of time-varying profits, which aligns with scenarios where survival rates decline over time. \citet{yu_team_2022} examine the OP with profits dependent on arrival and service times. They propose an exact algorithm based on Benders branch-and-cut, as well as a hybrid heuristic method, to solve the resulting nonlinear and nonconvex problem. A robust variant addressing random service times has also been explored \citep{yu_robust_2022}. In addition, the application of online learning to the OP is a growing area of research \citep{shiri_capacitated_2024,shiri_online_2020}. These studies assume no prior information about uncertain parameters, such as exact values, probability distributions, or uncertainty sets. Instead, parameters are incrementally revealed as vehicles visit locations, aligning well with the dynamic and uncertain nature of real-world disaster scenarios.

Other research focuses specifically on SAR planning. \citet{zheng_collaborative_2019} propose a collaborative human-UAV search planning framework aimed at minimizing the expected time for human rescuers to reach targets. \citet{liu_decision_2020} redefine rescue efficiency by introducing an intuitive metric to evaluate SAR operations. \citet{rauchecker_exact_2019} addressed the assignment and scheduling of rescue units to incidents under strict time constraints. They model this challenge as a generalized parallel machine scheduling problem, presenting an MIP formulation and a branch-and-price algorithm to generate efficient solutions.

\subsection{Machine Learning in SAR}

Machine learning methods, particularly RL, demonstrate significant potential in addressing complex real-world SAR planning problems, especially given their flexibility and rapid execution capabilities. Several studies leverage pointer networks \citep{bello2016neural} to tackle routing challenges. \citet{gama_reinforcement_2021} utilized RL to train pointer network models for solving the OP with time windows. \citet{wang_solving_2024} propose a double-layer optimization framework to solve large-scale OP. The external optimizer is based on a diversity evolutionary algorithm, and a pointer network is trained by deep RL as the inner optimizer. RL has also been applied to address complex and dynamic environments in SAR planning. \citet{lee_multi-start_2024} investigated the multi-start OP using deep RL, emphasizing the importance of near-optimal mission re-planning in real time to handle unexpected environmental changes. \citet{wan_deep_2024} propose an attention-based deep RL method for multi-UAV scheduling that accounts for time-varying profit. Their approach facilitates real-time decision-making in disaster data collection.

In summary, vehicle routing models along with machine scheduling models are widely used to formulate SAR tasks such as UAV search and routing, collaborative search operations, and SAR dispatch planning. Despite the computational challenges, these models effectively address time sensitivity and resource constraints, improving the efficiency and reliability of SAR operations. However, in scenarios requiring real-time decision-making and adaptive response strategies for dynamic disaster situations, RL has emerged as a complementary approach. While these data-driven techniques still face challenges related to interpretability, stability, and generalization, RL serves as a valuable addition to traditional optimization models in addressing the complexities of SAR tasks.

\section{Review of Relief Distribution} \label{LB:Relief}

Relief distribution is a critical component of disaster management, involving decisions on the allocation and transportation of resources under uncertain and dynamic conditions. These decisions aim to efficiently route vehicles, prioritize resources to meet urgent needs, and address challenges such as fluctuating demand patterns, disrupted infrastructure, and time-sensitive deliveries. State-of-the-art strategies frequently integrate real-time data and multi-modal transportation options to enhance the responsiveness and effectiveness of relief efforts.

\subsection{Optimization Model in Relief Distribution}

The relief distribution problem is frequently formulated as a vehicle routing problem (VRP). \cite{zhang2021humanitarian} proposed a combined drone OP for post-disaster logistics, where a truck and drone work together to assess disaster-affected areas. In this approach, the drone collects data from inaccessible or damaged locations, while the truck provides logistical support, including drone recharging and transportation between regions.  \cite{espejo2021multiagent} proposed a VRP formulation for humanitarian aid delivery, focusing on determining the optimal quantity of aid, vehicle routes, and delivery schedules from a central depot to multiple distribution points. 

Building on the foundational formulation like the VRP in relief distribution, more integrated models have emerged to address the complexities of real-world disaster scenarios. For instance, \cite{bayram2024joint} proposed a framework that integrates relief distribution with dynamic resource allocation, evacuation or shelter-in-place decisions, and staging strategies. Typically, this study simultaneously optimizes whether each zone should evacuate or shelter in place, allocates resources to road segments, and determines vehicle flows on road segments while ensuring the supply meets demand in a timely and efficient manner. \cite{li2020scenario} presented an integrated framework for hurricane relief logistics and casualty distribution, using multistage stochastic programming to optimize resource pre-positioning and evacuation flows. \citet{bhattaraimultistage} propose a multistage programming model for hurricane relief logistics, focusing on pre-positioning relief supplies, shelter locations, and determining distribution flows. Their approach aligns logistics operations with evolving demand to enhance efficiency and reduce costs.

Another important dimension of relief distribution is managing uncertainties, which include variability in demand for relief, disruptions to transportation networks, fluctuating resource availability, and potential facility failures \citep{desi2022optimization,wang2024two}. For example, \citet{siddig2024multistage} addressed evacuation-related uncertainties in disaster relief planning with a robust optimization model. They incorporate demand uncertainties, such as evacuee compliance and disaster characteristics, through event-based and box uncertainty representations. \citet{dalal2021robust} handle uncertainties in hurricane disaster relief logistics through a multistage stochastic programming model. They model evolving hurricane attributes with a Markov chain and adapt resource allocation across random planning stages.  \cite{balcik2020robust} address post-disaster routing with a robust optimization model that accounts for travel time variability. This approach provides feasible routing solutions under a wide range of possible transportation scenarios.

Equity in relief distribution is a critical challenge that demands careful consideration in disaster response. Ensuring equitable access to relief resources can prevent disparities in allocation. \cite{avishan2023humanitarian} proposed an adjustable robust optimization model that determines optimal routes, service times, and the sequence of site visits during disaster response. Their model emphasizes equitable distribution using a utility-based approach, where each site’s utility reflects the benefit from time and resources allocated. By balancing these utilities, the model ensures critical needs are met across all regions within the operational time constraints, avoiding resource concentration at easily accessible sites. \cite{sakiani2020inventory} proposed a multi-period inventory routing model for relief distribution, which integrates vehicle routing and network flow optimization. Equity is ensured by minimizing deprivation costs, defined as weighted shortages based on population size and area criticality, highlighting fairness in disaster relief.

\subsection{Machine Learning in Relief Distribution}

Machine learning emphasizes real-time adaptability in relief resource distribution, effectively managing dynamic demand, uncertainty, and equity. \cite{yu2021reinforcement} introduced a Q-learning framework modeling resource allocation as an MDP, where states represent allocation configurations, actions adjust resource distribution, and rewards optimize cost reductions. This approach efficiently balances delivery costs and deprivation penalties, demonstrating superior adaptability compared to traditional methods. Building on this foundation, \cite{fan2022dhl} employed a deep Q-network to address supply distribution in disaster scenarios with varying impacts. Deep Q-network incorporates real-time state-action mapping, enabling efficient allocation strategies that outperform heuristic and exact algorithms, particularly in large-scale, dynamic environments. Extending RL to real-time decision-making, \cite{wang2023deep} applied the deep deterministic policy gradient algorithm to rescue resource allocation during storm surges. In a case study of the 2018 Mangkhut storm surge, a deep deterministic policy gradient achieves near-optimal solutions with drastically reduced computation times, highlighting its suitability for adaptive logistics under tight constraints. \cite{van2023reinforcement} explored RL for last-mile relief distribution, integrating trucks and UAVs to mitigate travel time uncertainties caused by damaged infrastructure. The study optimizes delivery coverage and lateness penalties, showing significant performance improvements with UAV deployment, particularly in remote or low-demand areas.

\subsection{Simulation in Relief Distribution}


\cite{espejo2021multiagent} integrated ABM with an MIP model to optimize humanitarian aid distribution, accounting for behavioral factors such as victim impatience and interactions between victims and distribution points. A case study of the 2017 Mocoa landslide in Colombia demonstrates significant improvements in service levels and response efficiency through adaptive updates and dynamic optimization. Similarly, \cite{khalili2022stochastic} utilized simulation to estimate demand for essential relief supplies during cascading disasters. By analyzing interactions among vulnerable urban infrastructures, the model provides real-time inputs to a bi-level stochastic optimization framework, enabling precise supply allocation, minimizing costs, and maximizing coverage.


In essence, the field of relief distribution research focuses on the intricate logistics of delivering resources to disaster-affected areas, emphasizing efficiency, equity, and adaptability. Integrated optimization models are widely utilized to address the pre-positioning of supplies, evacuation planning, and routing under dynamic disaster conditions. Additionally, robust and stochastic optimization approaches have been developed to effectively manage uncertainties related to demand, transportation disruptions, and resource availability. Besides, machine learning is pivotal in enhancing real-time adaptability in relief distribution by effectively managing dynamic relief demand and damaged infrastructure. Furthermore, simulation techniques are employed to replicate real-world disaster scenarios, enabling the evaluation and refinement of relief strategies under various conditions, such as demand fluctuations, infrastructure failures, and behavioral factors.

\section{Conclusion and Future Research} \label{LB:conclusion}

This paper provides a comprehensive methodological review of five critical topics in disaster emergency response planning, focusing on optimization models, machine learning approaches, and simulation techniques. The inherent complexity of disaster response planning stems from the high uncertainty of disaster events, dynamic environmental conditions, and the need to incorporate human behavior. These factors make disaster planning an interdisciplinary challenge requiring innovative solutions that balance analytical rigor, computational efficiency, and practical adaptability.

Optimization models have traditionally served as the cornerstone of disaster emergency planning. Many problems in disaster emergency planning are formulated as classical optimization problems, including location-allocation, vehicle routing, optimal control, and scheduling. Optimization models have the advantage of delivering precise, often optimal, solutions to problems with well-defined objectives and constraints. However, To ensure tractability, these models often simplify the complexity of disaster scenarios by quantifying uncertainty through scenarios and relying on strong assumptions about human behavior. Additionally, since these problems are inherently NP-complete, the computational intractability of large-scale problems also limits their practical applicability. To address these limitations, machine learning and simulation methods have emerged as complementary approaches, offering improved realism in modeling disaster scenarios.

Machine learning offers an adaptive and data-driven approach to disaster response. Techniques such as supervised learning and RL enable systems to learn from historical data and improve decision-making in uncertain and dynamic environments. They are particularly effective in tasks requiring real-time adaptability compared with optimization models, such as predicting demand surges and optimizing casualty transport. However, adopting machine learning in disaster response faces challenges related to the availability and quality of training data, the interpretability of learned models, and the generalizability of solutions across different disaster scenarios. On the other hand, simulation techniques have emerged as powerful tools for capturing the complex nature of disaster environments. Methods such as ABM and discrete-event simulation allow researchers to model realistic scenarios that include interactions among evacuees, infrastructure disruptions, and cascading failures. Simulation provides critical insights into the emergent behavior of systems under stress and enables the evaluation of various planning strategies under different scenarios. However, the lack of analytical rigor in simulation approaches raises concerns about the weakness of guiding decision-making. Furthermore, simulations are computationally intensive, particularly when combined with optimization frameworks.

Based on our review, we identify several critical future research directions: (1) Rapid response to environmental changes: as the original operational strategy progresses, newly updated information requires corresponding planning adjustments. Online learning and RL demonstrate significant potential in this direction. (2) Accurate prediction in complex environments: Simulation and supervised learning are advantageous in this regard, but they primarily serve as evaluation tools. To optimize decision-making, these methods need to be integrated with complementary optimization frameworks. (3) Equity in disaster planning: Disaster planning is inherently human-centered, so ensuring equity by addressing the diverse needs of different individuals is a critical requirement. Multi-objective optimization can play a key role in balancing conflicting priorities to achieve equitable outcomes. (4) Quantifying and mitigating risk: Risk aversion is essential in disaster planning, as the potential impacts of extremely low-probability events cannot be overlooked. Probabilistic risk analysis and stochastic optimization methods can be applied to quantify risks and develop robust mitigation strategies. These directions highlight critical pathways for advancing disaster emergency response planning.

\appendix

\bibliographystyle{elsarticle-harv} 
\bibliography{cas-refs}





\end{document}